\documentclass[letterpaper, preprint, paper,11pt]{AAS}	

\usepackage{bm}
\usepackage{algorithm}
\usepackage{algorithmic}
\usepackage{amsmath}
\usepackage{amssymb}
\usepackage{subfigure}
\usepackage{float}
\usepackage[colorlinks=true, pdfstartview=FitV, linkcolor=black, citecolor= black, urlcolor= black]{hyperref}
\usepackage{overcite}
\usepackage{footnpag}			      	

\PaperNumber{24-199}

\begin{document}

\title{CONCURRENT OPTIMIZATION OF SATELLITE PHASING AND
TASKING FOR CISLUNAR SPACE SITUATIONAL AWARENESS}

\author{Malav Patel\thanks{PhD Candidate, Aerospace Engineering, Georgia Institute of Technology},  
Kento Tomita\thanks{PhD Candidate, Aerospace Engineering, Georgia Institute of Technology}
\ and Koki Ho\thanks{Dutton-Ducoffe Professor and Associate Professor, Aerospace Engineering, Georgia Institute of Technology}
}

\maketitle{}

\begin{abstract}
Recently, renewed interest in cislunar space spurred by private and public orga-
nizations has driven research for future infrastructure in the region. As Earth-
Moon traffic increases amidst a growing space economy, monitoring architectures
supporting this traffic must also develop. These are likely to be realized as con-
stellations of patrol satellites surveying traffic between the Earth and the Moon.
This work investigates the concurrent optimization of patrol satellite phasing and
tasking to provide information-maximal coverage of traffic in periodic orbits.
\end{abstract}

\section{Introduction}
Constellation design and optimization is of fundamental importance in space systems engineering, driven primarily by cost and maintenance. The design problem for Earth-based constellations and a variety of objectives is well understood in the literature\cite{luo2017constellation,lee2020satellite, ulybyshev2008satellite, shimane2024orbital, guan2020optimal, davis2012reducing}. Recently, a growing interest in missions beyond low Earth orbit has prompted a flurry of research into cislunar space, defined as the space between the Earth and Moon. Because of the vast distances involved, conventional infrastructure used in space systems near the Earth do not translate well to the cislunar regime\cite{frueh2021cislunar}.

One current research thrust investigates how to achieve space situational awareness. This involves (but is not limited to) monitoring traffic moving within the region, transfers between lunar and Earth-based orbits, and establishing robust communication architectures.
Space situational awareness in cislunar space can be realized with ground-based sensors, space-based sensors, or a combination of both, where space-based sensors reside on suitable orbits\cite{gupta2021earth}. Recent work has begun to investigate a fusion of space and ground based surveillance networks to offset the weaknesses of exclusive ground based observation networks\cite{harris2021expanding}. Ground-based sensors struggle to observe this vast region due to cloud cover, restricted field of view, and insufficient observational power\cite{2021spde.confE.188F}.  Space-based sensors are free from the first two drawbacks and have competitive observational power. In this work, we investigate constellations of space-based observers. 

One facet of space situational awareness is the ability to monitor targets in cislunar space on periodic orbits. We represent each target as a tuple ($\boldsymbol{x}(t), \boldsymbol{P}(t)$) of its state and covariance. The time evolution of the state and covariance are governed by the dynamics of the underlying gravitational model and the observations provided by space-based observers, which is provided by common filtering algorithms. This work uses the circular restricted 3 body problem (CR3BP) as the gravitational model and the extended kalman filter (EKF) for the filtering algorithm. 

Many works optimize either the phasing of satellites in a constellation for an objective while ignoring how to task them or vice versa: assuming a fixed constellation and optimizing for the best tasking schedule. Hartigan et al.\cite{hartigan2023optimization} optimize the phase of a small constellation of satellites near the lunar south pole to provide necessary foundational infrastructure for the development of larger constellations. Vendl et al.\cite{vendl2021cislunar} present a method to maximize observation capability of a much larger cone of points in cislunar space by optimizing the phase of an observer. Visonneau et al. \cite{Visonneau_2023} present a hidden genes genetic algorithm to optimize a set of observers and their phases in the CR3BP, aiming to minimize the total number of observers, the stability index of chosen orbits, and maximize the total coverage of a region of cislunar space. Patel et al.\cite{patel2023cislunar} investigate how to optimally phase a set of observers on a catalog of orbits to satisfy spatio-temporal demand. However, these works do not consider how to task each observer, only optimizing a metric that indicates an observer is \emph{available} to monitor a target at a specified time. On the other side, Fedeler et al. \cite{FEDELER2022792, FEDELER20245266} illustrate methods for cooperative and decentralized tasking of space based observers for tracking space objects. However, they keep the number of observers and their phases fixed.  It is not unreasonable to claim that there may exist some coupling between the phasing of the satellites and the resulting optimal tasking. Fahrner et al.\cite{fahrner2022capacity} present a method to optimize capacity (sensor scheduling) and observability (constellation design) simultaneously using a genetic and greedy algorithms. This work aims to concurrently optimize satellite phasing and the tasking strategy using a two-level optimization scheme that relies on linear programming and gradient-based optimization.

The work is presented as follows. First we present some background on the gravitational model used. Next we present the method for concurrent optimization of a constellation's phase and tasking of its sensors, presenting two formulations with different objectives. Then we provide greedy and exhaustive algorithms for optimization. We show empirically that gradient based optimization is well-suited for part of the solution process. Finally, we present 3 case studies. The first compares the performance of the aforementioned algorithms. The second compares the performance of myopic and optimal policies. The third compares the results from using the two aforementioned objectives.

\section{Background}
\subsection{Circular Restricted 3 Body Problem}
This work requires a gravitational model to propagate orbits for agents and targets. Here we use the circular restricted 3 body problem (CR3BP) as a model for the dynamics, since it offers a reasonable approximation to the motion of a small satellite influenced by the gravity of two larger masses. The equations of motion are given by 

\begin{equation}
\begin{aligned}
    \ddot{x} - 2\dot{y} &= \frac{dU}{dx}\\
    \ddot{y} + 2\dot{x} &= \frac{dU}{dy}\\
    \ddot{z} &= \frac{dU}{dz}
\end{aligned}
\end{equation}

where $U$ is an effective potential driven by the rotating frame of the two larger masses, given by
\begin{equation}
    U = \frac{x^2 + y^2}{2} + \frac{1-\mu}{r_1} +\frac{\mu}{r_2}
\end{equation}

Here, $\mu = \frac{m_2}{m_1 + m_2}$ is the mass ratio of the two larger bodies with masses of $m_1$ and $m_2$, and $r_1$ and $r_2$ denote the distances of the satellite to each body.

\begin{equation}
    \begin{aligned}
        r_1 &= \sqrt{(x+\mu)^2 + y^2 + z^2}\\
        r_2 &= \sqrt{(x-(1-\mu))^2 + y^2 + z^2}
    \end{aligned}
\end{equation}

Initial conditions are propagated numerically using the Heyoka library\cite{Biscani_2021}.

\section{Method}
\subsection{Quantifying Space Situational Awareness}
The ability to monitor or provide space situational ``awareness" can be quantified using information theory. More specifically, observations by satellites provide information and tend to reduce the eigenvalues of the covariance matrix $\boldsymbol{P}(t)$. Tomita et al.\cite{tomita2023multispacecraft} provide an information-theoretic approach to satellite constellation design in order to monitor targets on periodic orbits in the CR3BP. They formulate 2 integer linear programs to solve for the control necessary to maximize either (1) the cumulative information gain or (2) the minimum target information gain over a given simulation time. For brevity, the programs are outlined below.

\subsubsection{Maximizing Cumulative Information (Max Objective)}

\begin{equation}
\begin{aligned}
   \max_{u_{ijk}}& \sum_{i=1}^{M} \sum_{j=1}^{N} \sum_{k=0}^{L-1} u_{ijk}\  \text{tr} \big(I_{ij}(t_{L}, t_{k})\big), && \quad I_{ij}(t_{L}, t_{k}) = \Phi_{j}(t_{k}, t_{L})^{\top} H_{ijk}^{\top} R_{ij}^{-1} H_{ijk} \Phi_{j}(t_{k}, t_{L}) \\\
   \text{s.t.}& \quad \sum_{j=1}^{N}u_{ijk} \leq 1\quad \forall i, k && \text{(observer cannot view more than 1 target at each timestep)} \\\
   & \quad u_{ijk} \in \{0, 1\}
\end{aligned} \label{eq:1}
\end{equation} 

\subsubsection{Maximizing Minimum Target Information (MaxMin Objective)}

\begin{equation}
\begin{aligned}
   \max_{u_{ijk}} \min_{j}& \sum_{i=1}^{M} \sum_{k=0}^{L-1} u_{ijk}\  \text{tr} \big(I_{ij}(t_{L}, t_{k})\big), && \quad I_{ij}(t_{L}, t_{k}) = \Phi_{j}(t_{k}, t_{L})^{\top} H_{ijk}^{\top} R_{ij}^{-1} H_{ijk} \Phi_{j}(t_{k}, t_{L}) \\\
   \text{s.t.}& \quad \sum_{j=1}^{N}u_{ijk} \leq 1\quad \forall i, k && \text{(observer cannot view more than 1 target at each timestep)} \\\
   & \quad u_{ijk} \in \{0, 1\}
\end{aligned} \label{eq:2}
\end{equation}

Here, $i$ indexes the observers, $j$ indexes the targets, and $k$ indexes the timesteps. $ I_{ij}(t_{L}, t_{k})$ is the information matrix for target $j$ observed by observer $i$ at timestep $t_k$ propagated to the final time $t_L$. This formulation solves the control problem conditioned on a given set of observers. In other words, given the phase of the observers, we can retrieve the optimal control that maximizes the information gain. However, we may also consider the initial phase of each observer as a design variable and choose them such that the information gain is maximized.

\subsection{Concurrent Tasking and Phasing Problem}
Considering both the phase and the control as design variables in a concurrent optimization problem, we present the following program when the objective is to maximize the cumulative information.

\begin{equation}
\begin{aligned}
    \underset{\boldsymbol{x}}{\max}& &&\sum_{i=1}^{M} \sum_{j=1}^{N} \sum_{k=0}^{L-1} u^*_{ijk}\  \text{tr} \big(I_{ij}(t_{L}, t_{k}, \boldsymbol{x})\big) \\\
   \text{s.t.}& \quad &&0 \leq \boldsymbol{x} \leq 1 \\
   &&&\quad u^*_{ijk} \in
   \begin{aligned}[t] \underset{u_{ijk}}{\arg\max}& &&\sum_{i=1}^{M} \sum_{j=1}^{N} \sum_{k=0}^{L-1} u_{ijk}\  \text{tr} \big(I_{ij}(t_{L}, t_{k}, \boldsymbol{x})\big) \\\
   \text{s.t.}& \quad &&\sum_{j=1}^{N}u_{ijk} \leq 1\quad \forall i, k  \\\
    &&&\quad u_{ijk} \in \{0, 1\}
\end{aligned} \\\
\end{aligned} \label{eq:3}
\end{equation}

Here, $\boldsymbol{x}[i] \in [0, 1)$ captures the initial phase of the $i$-th observer. $\boldsymbol{x}[i]$ enters the objective function through the observation jacobian $H_{ijk}$. This jacobian is dependent on the position of the observer at time step $k$, and the observer position at timestep $k$ is in turn dependent on initial phase $\boldsymbol{x}[i]$. 

This program consists of an upper level which solves for the optimal phasing of observers $\boldsymbol{x}^*$, and a lower level which solves for the optimal control $u^*_{ijk}$ given a choice of $\boldsymbol{x}$. If we instead would like to consider the minimum target information as the objective, we simply replace the upper-level and lower-level objectives accordingly. We hypothesize the maxmin objective can help diversify the allocation of sensor resources amongst more targets when compared to the max objective. A case study is presented later in the Results section. We note that of the two algorithms presented in this work, the greedy algorithm returns the optimal solution only when we use the max objective (i.e. Equation \ref{eq:1}). To guarantee optimality with the maxmin objective, we must opt for an algorithm that searches the entire parameter space.
\subsection{Greedy Solution Approach}
The program in ($\ref{eq:3}$) can be solved efficiently using a greedy algorithm that optimizes the phase of each observer independently. To see this we first note the matrix $I_{ij}(t_{L}, t_{k})$ has non-negative coefficients, $\boldsymbol{A}[i,j,k] \triangleq \text{tr}\big(I_{ij}(t_{L}, t_{k})\big) \geq 0$. So, the solution to the lower level program \textit{without} the constraint $\sum_{j=1}^{N}u_{ijk} \leq 1\ \  \forall i, k$ is trivially a 3d tensor of all ones, $\boldsymbol{u}^* = \boldsymbol{1}$. The constraint forces the array $\boldsymbol{u}^*[i, :, k]$ to be a one-hot vector $\forall \ \ i, k$. Lastly, notice that $\boldsymbol{A}[i, :, k]$ depends on the position of observer $i$ at timestep $k$ and not on the position of any other observer. As a result, $\boldsymbol{A}[i, :, :]$ and $\boldsymbol{u}^*[i, :, :]$ are functions of $\boldsymbol{x}[i]$ only\footnote{There is a dependence on the target positions, but this does impact the results of the discussion.}. Define $f(\boldsymbol{x}[i]; \boldsymbol{u}^*) \triangleq \boldsymbol{u}^*[i, :, :]\  \cdot \boldsymbol{A}[i, :, :]$, where ``$\cdot$" represents the dot product. With this information, we can rewrite the program in (\ref{eq:3}) as
\begin{equation}
\begin{aligned}
    \underset{\boldsymbol{x}}{\max}& &&\sum_{i=1}^{M} f(\boldsymbol{x}[i];\boldsymbol{u}^*) \\\
   \text{s.t.}& \quad &&0 \leq \boldsymbol{x} \leq 1 \\
   &&&\quad \boldsymbol{u}^* \in
   \begin{aligned}[t] \underset{\boldsymbol{u}}{\arg\max}& && \boldsymbol{u} \cdot \boldsymbol{A}, \quad \boldsymbol{A}[i,j,k] \triangleq \text{tr}\big(I_{ij}(t_{L}, t_{k}, \boldsymbol{x})\big) \\\
   \text{s.t.}& \quad &&\sum_{j=1}^{N}\boldsymbol{u}[i, j, k] \leq 1\quad \forall i, k  \\\
    &&&\quad \boldsymbol{u}[i, j, k] \in \{0, 1\}
\end{aligned} \\\
\end{aligned} \label{eq:4}
\end{equation}
Note that the upper and lower level objectives are identical functions, but written in more convenient forms. In this form, it is easy to see that maximizing the upper level objective is equivalent to maximizing each term in the summation independently. As a result, the optimal decision vector $\boldsymbol{x}^*$ can be constructed by sequentially optimizing each of its entries instead of searching a multidimensional space. The greedy method also reduces the size of the lower level linear program since we only consider one observer while we solve Equation \ref{eq:3}. With one observer being solved at a time, the number of binary variables in the lower level linear program is reduced by a factor of $M$, the number of total observers, compared to the the lower level linear program in the exhaustive approach. Note that this decomposition is not valid if we choose our objective to maximize the minimum target information. In this case, an optimizer must search the multidimensional space. Algorithms 1 and 2 below show pseudocode for the exhaustive and greedy methods.

\begin{minipage}[t]{0.45\textwidth}
    \begin{algorithm}[H]
    \caption{Exhaustive Search}
    \begin{algorithmic}[1]
        \STATE O $\gets$ Initialize Observer set
        \STATE SE $\gets$ Initialize space environment with targets and observers.
        \STATE $\boldsymbol{x}^* \gets$  solve (\ref{eq:3}) with M = $|\text{O}|$ via BFGS (upper level) and LP (lower level)

        \STATE Return $\boldsymbol{x}^*$
    \end{algorithmic}
    \end{algorithm}
\end{minipage}
\hfill
\begin{minipage}[t]{0.45\textwidth}
    \begin{algorithm}[H]
    \caption{Greedy Search}
    \begin{algorithmic}[1]
        \STATE O $\gets$ Initialize Observer set
        \STATE SE $\gets$ Initialize empty space environment with targets.
        \STATE $\boldsymbol{x}^* \gets$ []
        \FOR{o in O}
            \STATE SE.remove\_all\_observers()
            \STATE SE.add\_observer(o)
            \STATE $x_\text{o} \gets$ solve (\ref{eq:3}) with M = 1 via BFGS (upper level) and LP (lower level)
            \STATE $\boldsymbol{x}^* \gets$ [$\boldsymbol{x}^*\ \  x_\text{o}$]
        \ENDFOR
        \STATE Return $\boldsymbol{x}^*$
    \end{algorithmic}
    \end{algorithm}
\end{minipage}

\section{Results}
In this section, we present empirical results supporting the use of gradient-based optimization for the upper-level objective. An ablation against this method using derivative-free optimization is left for future work.


Then we show empirically why a greedy algorithm is suited for this optimization problem. We then show three case studies examining the performance of exhaustive/greedy algorithms, myopic/optimal policies, and max/maxmin objectives. Orbits for constellation design in these sections are provided by NASA JPL\cite{JPL_catalog} and work from Broucke. \cite{1968port.book.....B}

\subsection{Gradient Based Methods and Smooth Objectives}
\begin{figure}[H]
  \centering
  \begin{minipage}[b]{0.45\textwidth}
    \includegraphics[width=\textwidth]{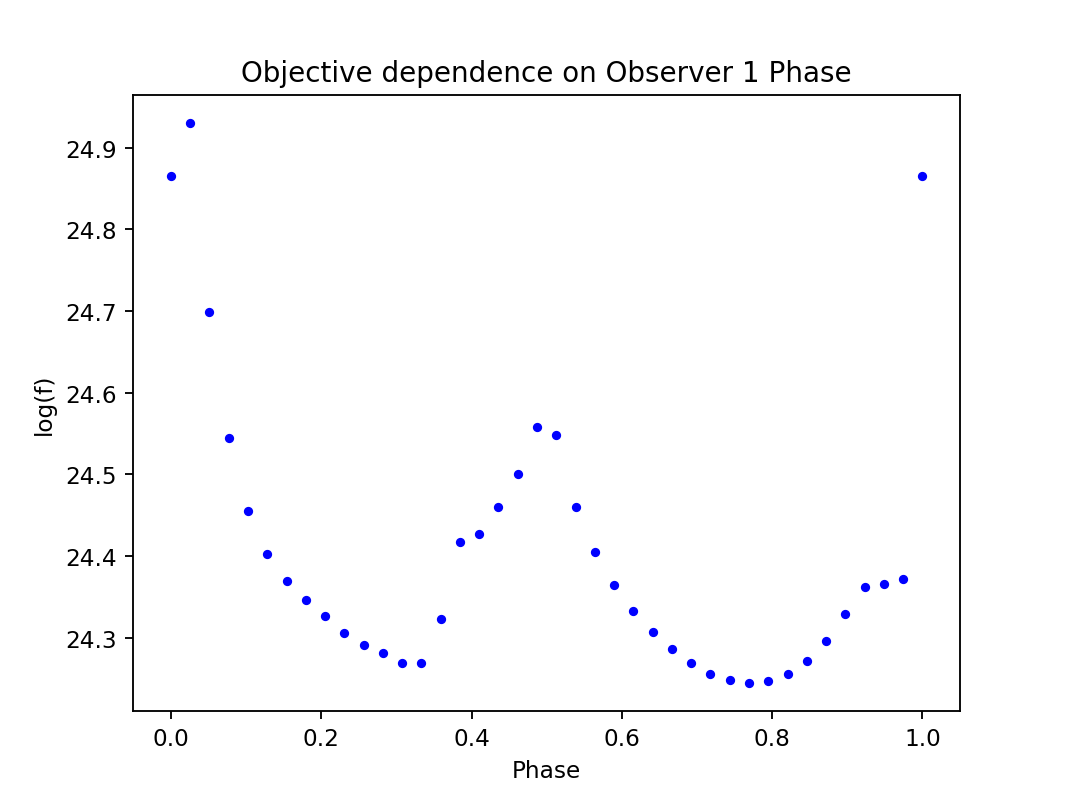}
  \end{minipage}
  \hfill
  \begin{minipage}[b]{0.45\textwidth}
    \includegraphics[width=\textwidth]{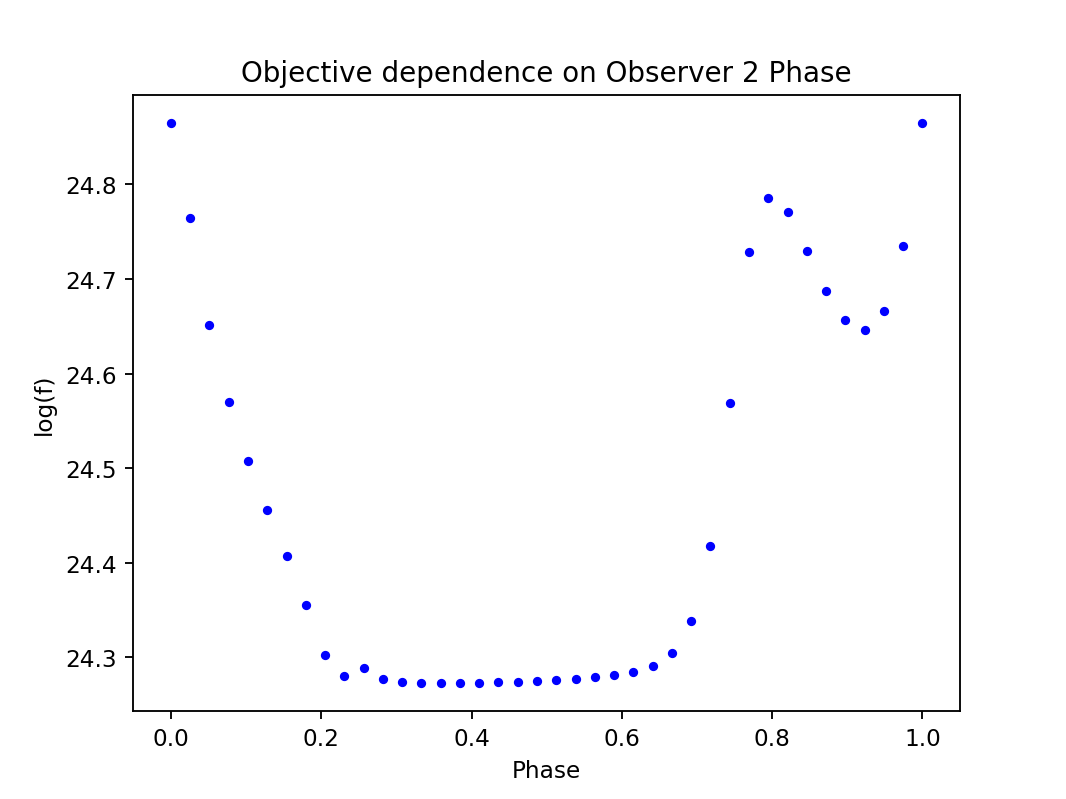}

  \end{minipage}
    \caption{The logarithm of the objective plotted against the phase of the first observer (left) and the second observer (right). The phase is parameterized as a fraction of the period of the observer.}
    \label{fig:1}
\end{figure}

Empirically we find the lower level program can be solved in seconds using commercial software\cite{gurobi} for problems with $\mathcal{O}(1)$ observers, $\mathcal{O}(1)$ targets, and $\mathcal{O}(10^2)$ timesteps. Furthermore, we find that the objective is well-behaved as the design vector $\boldsymbol{x}$ is varied. See Figure \ref{fig:1} for an illustration of how the objective varies with the design vector for a simple case of two observers tasked with monitoring two targets.

Fast lower level solve times suggest that the upper level objective can be evaluated inexpensively. Furthermore, continuity of the upper level objective with respect to $\boldsymbol{x}$ suggest that estimates of the gradient and curvature help accelerate convergence to a stationary point. As a result, techniques like gradient descent, Newton's method, and their relatives can accelerate convergence when compared to derivative free optimization algorithms. We use existing implementations of these algorithms in the PyGMO library\cite{Biscani2020}.

\begin{figure}[H]
    
  \centering
  \begin{minipage}[b]{0.49\textwidth}
    \includegraphics[width=\textwidth]{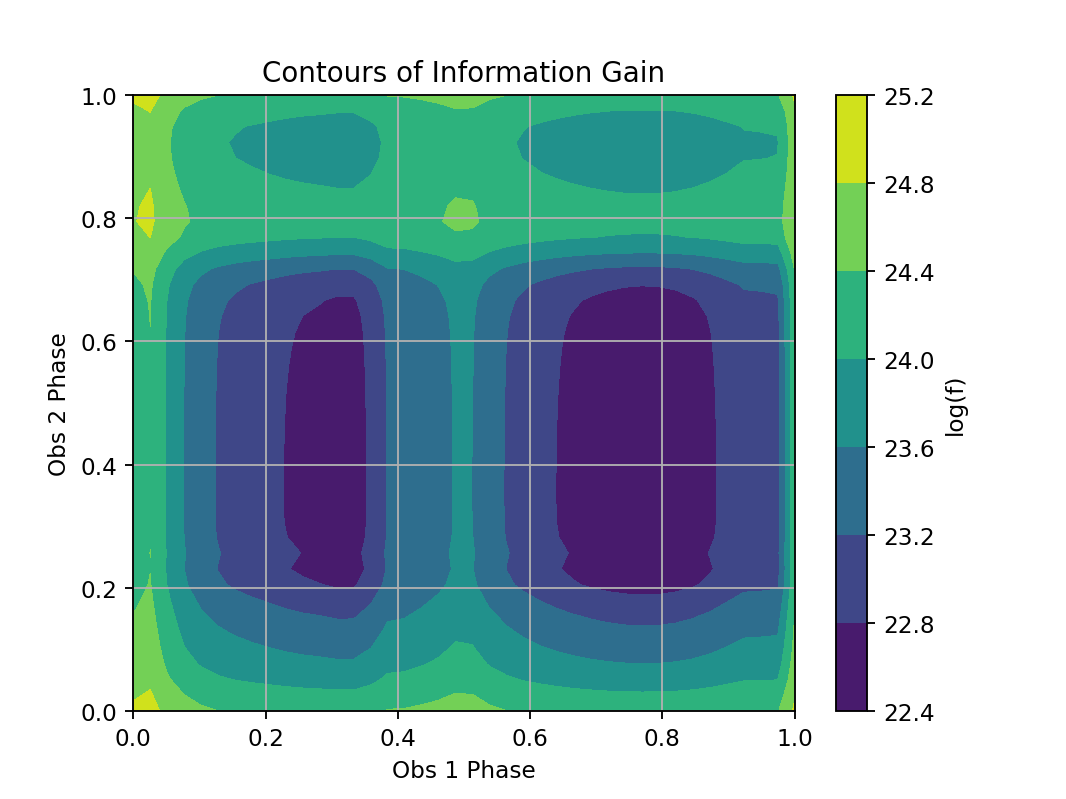}
  \end{minipage}
  \hfill
  \begin{minipage}[b]{0.49\textwidth}
    \includegraphics[width=\textwidth]{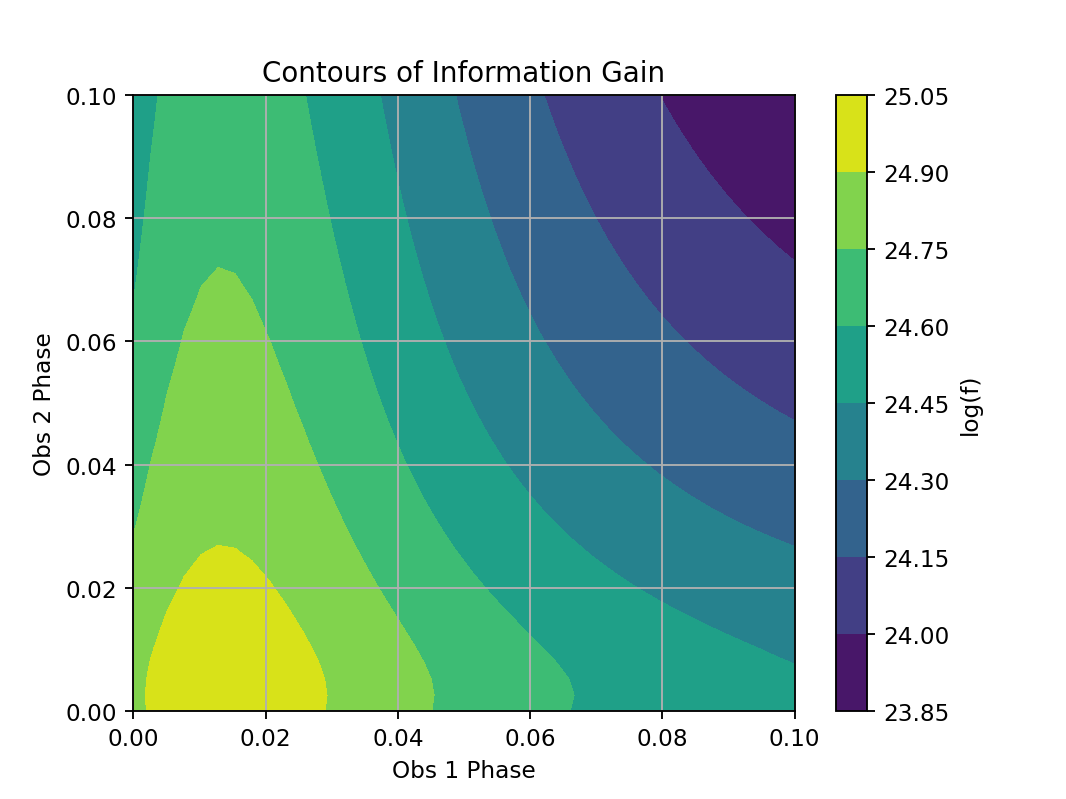}

  \end{minipage}
    \caption{Contours of the objective in the entire search space (left) and near the optimizer (right)}
    \label{fig:3}
\end{figure}

For a pair of observers and a small number of targets (i.e. the example in Figure \ref{fig:1}), the L-BFGS method is fast, converging to the optimal phase $\boldsymbol{x}^* = [0.01343, 0.0023]$ within 4 seconds. At this scale, we can verify the solution by looking at contours of the objective, shown in Figure \ref{fig:3}. However, the search space (and thus the algorithm run time) grows with the number of observers. At larger scales, a gradient based approach in an increasingly high dimensional space will prove costly. Fortunately, for the maximum cumulative information objective, the greedy decomposition reduces the search to a sequence of $N$ independent searches in the unit interval $[0, 1]$.
\subsection{Greedy Approach Validation}
\begin{figure}[H]
\centering\includegraphics[width=0.95\textwidth]{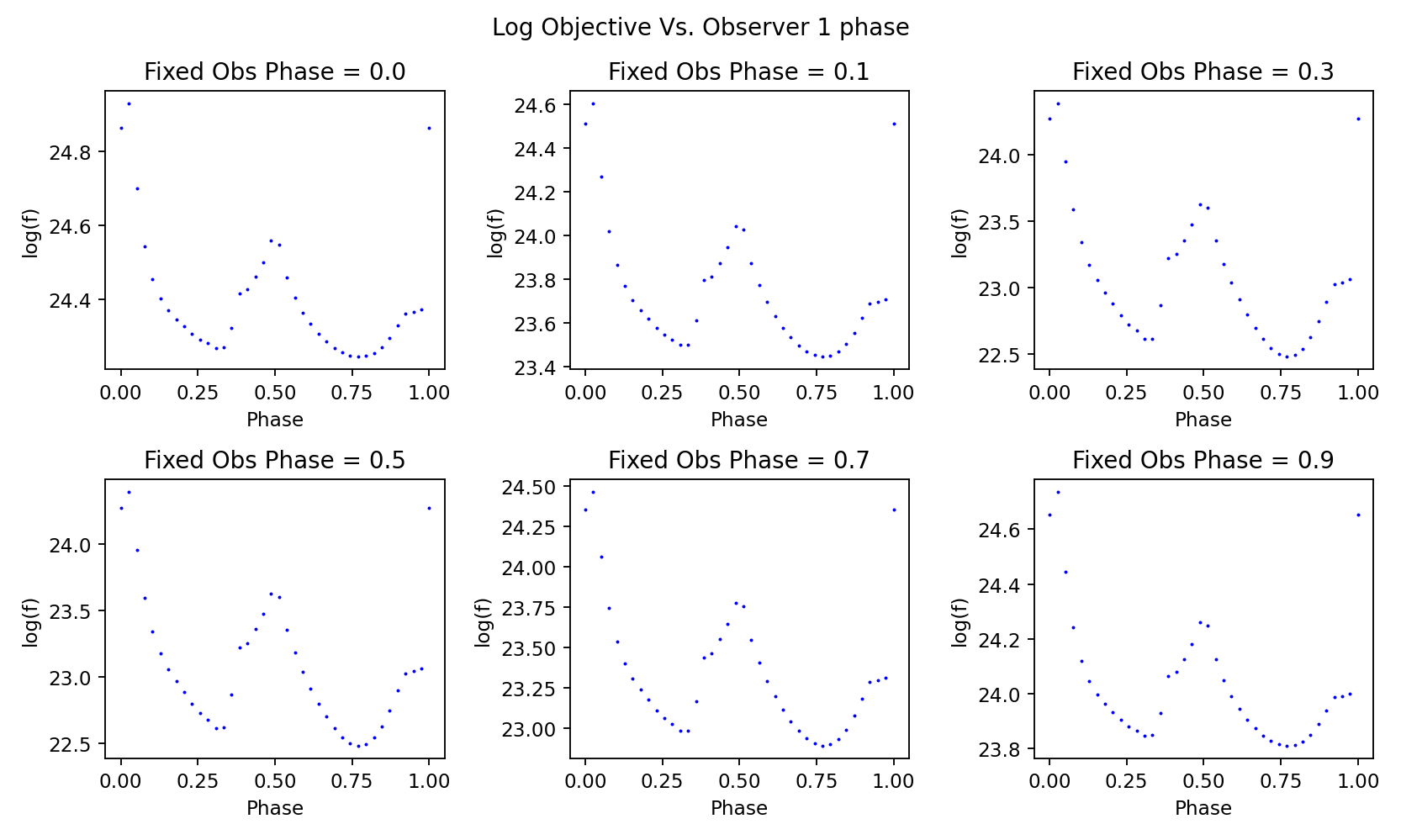}
	\caption{The logarithm of the objective plotted against the phase of an observer. The phase of the second observer is held constant. Although the scale of the information gain changes, the location of the optimal phase for observer 1 remains the same.}
	\label{fig:2}
\end{figure}

Empirically, the separability of the objective can be seen in Figure \ref{fig:2}, which plots the logarithm of the objective of one observer against its phase while the phase of the second observer is held constant. The characteristic "W" shape of the objective persists as the phase of the second observer is varied across each of the 6 plots. Although the scale of the objective changes, the location of the optimal phase for observer 1 remains independent of the phase of observer 2. 

\subsection{Case Study 1: Comparing Exhaustive and Greedy Approach}
The results have verified that (1) gradient based methods are well-suited for the upper level optimization problem and (2) the greedy decomposition in Equation \ref{eq:4} of the optimization problem in Equation \ref{eq:3} are equivalent problems. To compare the two approaches, we present an environment with 4 observers and 3 targets arranged on a variety of periodic orbits in the CR3BP. The environment is shown below in Figure \ref{fig:4}. 

Note that care was taken to choose orbits where no agent and target has an opportunity to collide, regardless of their phases. This was done primarily to avoid exploding coefficients in the lower level linear program defined in Equation \ref{eq:3}, since the observation jacobian $H$ contains several terms proportional to the reciprocal separation between the target and agent. We show in a later section that this effect may become dominant when there are close approaches between a target and observer, and as a result a myopic policy can perform close to optimally.

\begin{figure}[H]
  \centering
  \begin{minipage}[b]{0.6\textwidth}
    \includegraphics[width=\textwidth]{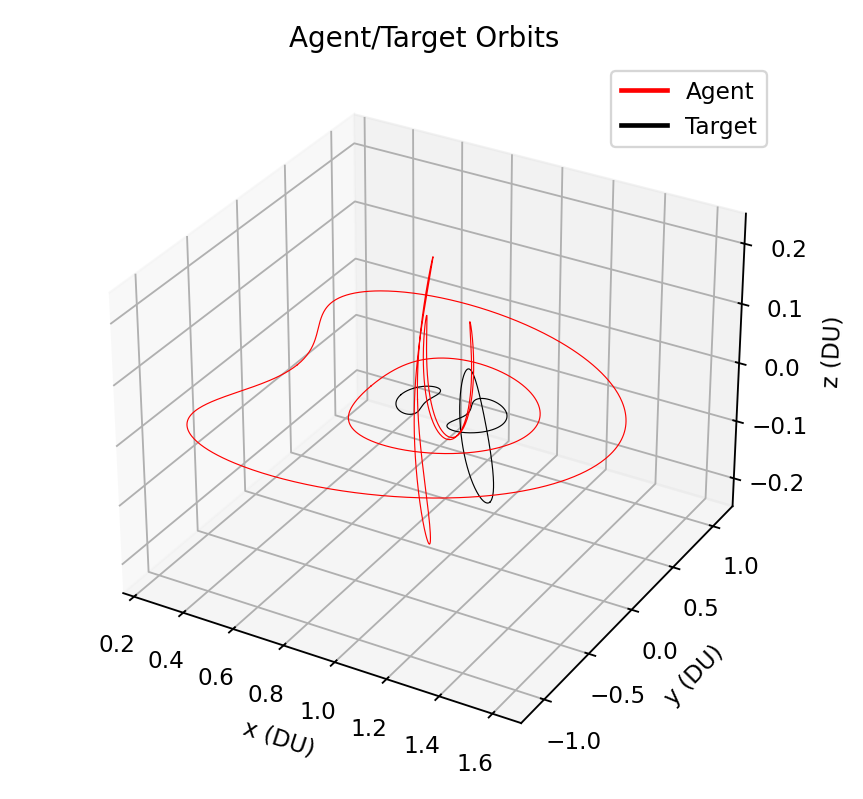}
  \end{minipage}
  \hfill
  \begin{minipage}[b]{0.39\textwidth}
    \includegraphics[width=\textwidth]{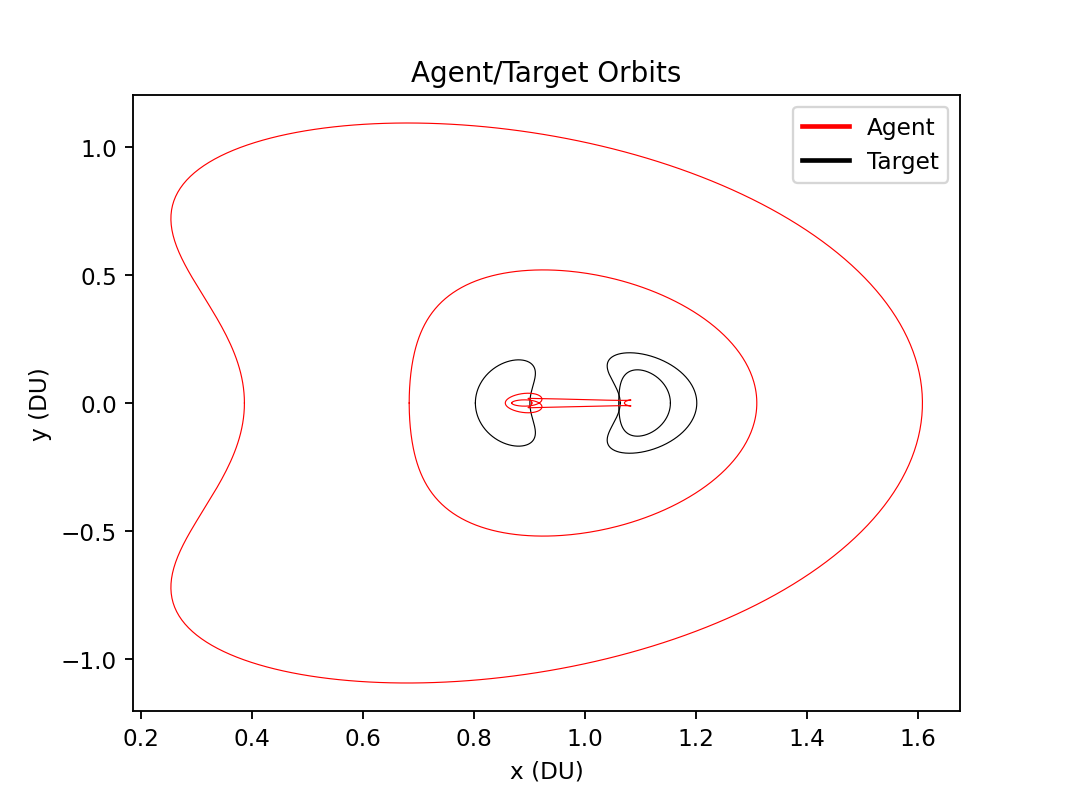}
    \includegraphics[width=\textwidth]{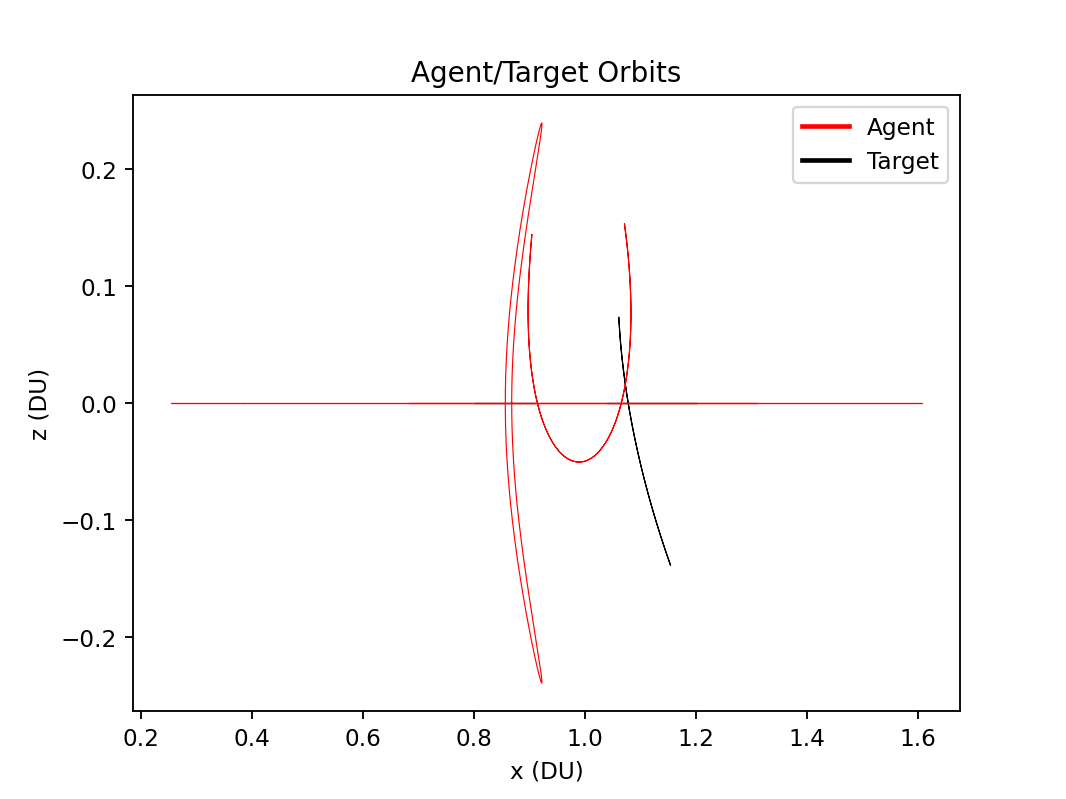}

  \end{minipage}
    \caption{Agent/Target orbits (left) and their projections on the xy plane (top right) and xz plane (bottom right). Orbits were chosen carefully to ensure that no intersection occurs between an agent and target orbit.}
    \label{fig:4}
\end{figure}

\subsubsection{Max Objective}
We choose an initial guess $\boldsymbol{x}_0 = [0.5, 0.5, 0.001, 0.7]$ motivated primarily by producing plots similar to those in Figure \ref{fig:1} and choosing a point close to a perceived maximizer. Table \ref{tab:1} shows the results of each approach. 

\begin{table}[htbp]
	\fontsize{10}{10}\selectfont
    \caption{Greedy vs. Exhaustive Search Results and Performance on Max Objective}
   \label{tab:label}
        \centering 
   \begin{tabular}{c  | c | c | c } 
      
                  & $\boldsymbol{x}^*$ & Solve Time (sec) & log(f)\\
      \hline 
      Greedy      & [0.4819 , 0.4764 , 0.0000 , 0.7064] & 5.07 & 33.9340714 \\
      Exhaustive   & [0.4794 , 0.4768 , 0.0003 , 0.7064] & 28.67 & 33.9340712 \\
      
   \end{tabular}
   \label{tab:1}
\end{table}


Practically speaking, the results of both approaches show that an optimizer exists in the neighborhood of $\boldsymbol{x} = [0.48, 0.48, 0.0, 0.706]$. However, we the greedy algorithm returns a solution with a log objective that is larger by 2e-7. Although the difference is negligible, it does change the optimal control returned by the two algorithms. The control differs at time steps 175 and 208. At time step $k = 175$, the greedy algorithm suggests tasking the second observer to the second target while the exhaustive algorithm suggests tasking that observer to the third target. At time step $k = 208$, the greedy algorithm suggests tasking the first observer to the second target while the exhaustive algorithm suggests tasking that observer to the third target. Numerically, the greedy algorithm returns a marginally better solution. For practical constellation design purposes, this advantage is negligible, but the greedy approach still boasts a 5x speedup in solve time over the exhaustive algorithm.


\subsubsection{MaxMin Objective} The greedy decomposition is no longer valid if we choose the maxmin objective because the $\min_j$ operation in Equation \ref{eq:2} introduces coupling between the decisions of agents (i.e. agents can no longer act independently to obtain the optimal solution). Table \ref{tab:2} shows the results of the two algorithms using the maxmin objective. We use an initial guess $\boldsymbol{x}_0 = [0.5, 0.5, 0.6, 0.2]$.

\begin{table}[htbp]
	\fontsize{10}{10}\selectfont
    \caption{Greedy vs. Exhaustive Search Results and Performance on MaxMin Objective}
   \label{tab:label}
        \centering 
   \begin{tabular}{c  | c | c | c } 
      
                  & $\boldsymbol{x}^*$ & Solve Time (sec) & log(f)\\
      \hline 
      Greedy      & [0.5444, 0.5214, 0.6723, 0.2193] & 54.95 & 22.9978 \\
      Exhaustive   & [0.4958, 0.5042, 0.6033, 0.2067] & 93.02 & 23.1169 \\
      
   \end{tabular}
   \label{tab:2}
\end{table}

We see a non-negligible gap in the objective between the algorithms, with the greedy algorithm performing sub-optimally, as expected. Furthermore, the returned control is vastly different between the two algorithms. Of the 215 time steps in each control array, 208 steps had different control.

In both max and maxmin formulations of the optimization problem, the objective landscape is non-convex, making it difficult to find a globally optimal solution. However, one can reduce the probability of convergence to local optima by executing numerous instances of gradient descent, each at a different initial condition. The choice of initial conditions can be driven heuristically by information from plots like Figure \ref{fig:1} or driven by DoE methods. An analysis of marginal improvement to the objective against the number of gradient descent threads is left to future work.

\subsection{Case Study 2: Sub-Optimality of the Myopic Policy}
In this section we show an example of how a myopic policy can be close to optimal when agent and target orbits are chosen to have intersections or close approaches to target orbits. However, the optimal policy is considerably better when this is not the case.

A myopic policy defines control for each agent at a particular timestep that is solely determined by the distances to each target. A myopic policy $u^{\text{myop}}$ would have the form,

\begin{align}
    u_{ijk}^{\text{myop}} = \min_{j}\  \{ \| \boldsymbol{r}_i(t_k) - \boldsymbol{r}_j(t_k) \| \}
\end{align}

\subsubsection{Intersecting Orbits} 
\begin{figure}[H]
    
  \centering
  \begin{minipage}[b]{0.49\textwidth}
    \includegraphics[width=\textwidth]{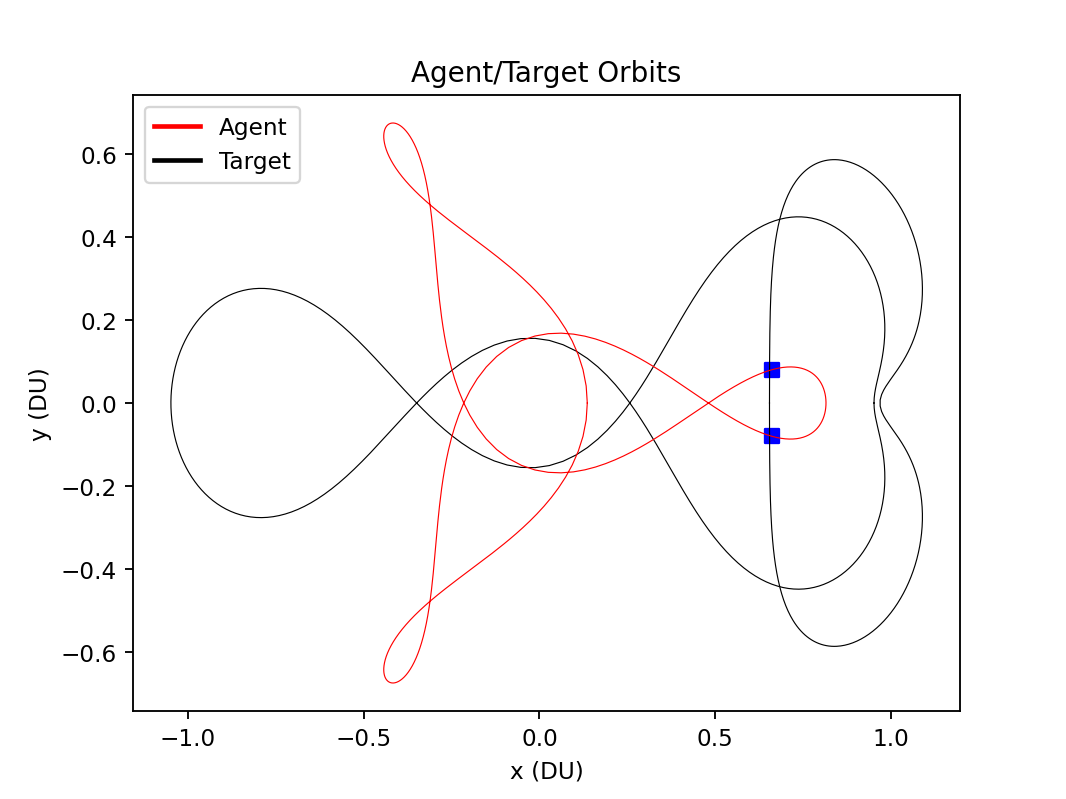}
  \end{minipage}
  \hfill
  \begin{minipage}[b]{0.49\textwidth}
    \includegraphics[width=\textwidth]{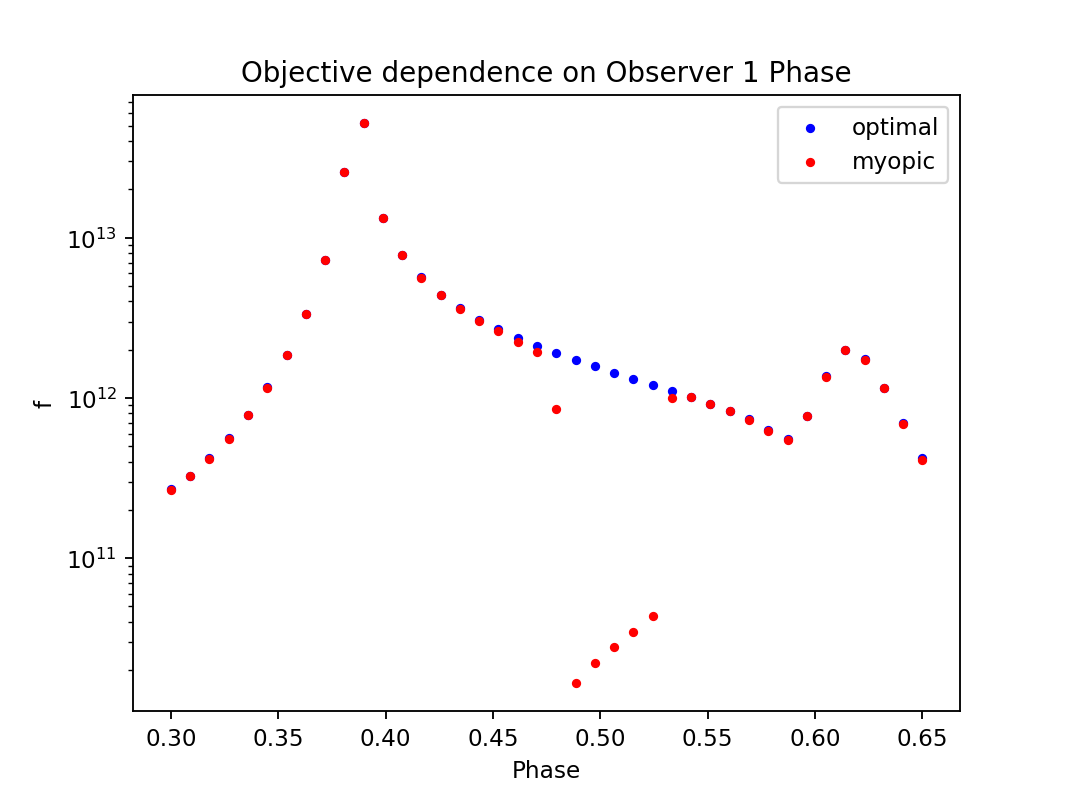}

  \end{minipage}
    \caption{Agent/Target Orbits (left) chosen intentionally to intersect. Objective dependence on the observer phase (right). Peaks correspond to square markers on the orbit plot (left), where an agent and target intersect.}
    \label{fig:6}
\end{figure}

Suppose that orbits are chosen such that at least one agent and target orbit intersect. Then the optimizer would unfavorably choose a solution $\{\boldsymbol{x}^*, \boldsymbol{u}\}$ that places agents and targets on collision courses. This is primarily because the observation jacobian $H$ has entries proportional to $\frac{1}{\| \boldsymbol{r}_T - \boldsymbol{r}_O \|^3}$, where $\boldsymbol{r}_T$ is the position of the target and $\boldsymbol{r}_O$ is the position of the agent. Hence, configurations leading to a small distance between agents and targets will explode the linear program's objective coefficients. Empirically we find that at the optimal phase $\boldsymbol{x}^*$, the objectives achieved by myopic and optimal control become practically identical. In this case, there is little justification for the extra leg-work required to determine the optimal control when a myopic policy comes within 1\% of the true objective. Figure \ref{fig:6} shows the myopic and optimal policy performance for a configuration of 1 agent and 2 targets with intersections between the orbits.

\subsubsection{Non-Close/Close Approaches}
\begin{figure}[H]
    
  \centering
  \begin{minipage}[b]{0.49\textwidth}
    \includegraphics[width=\textwidth]{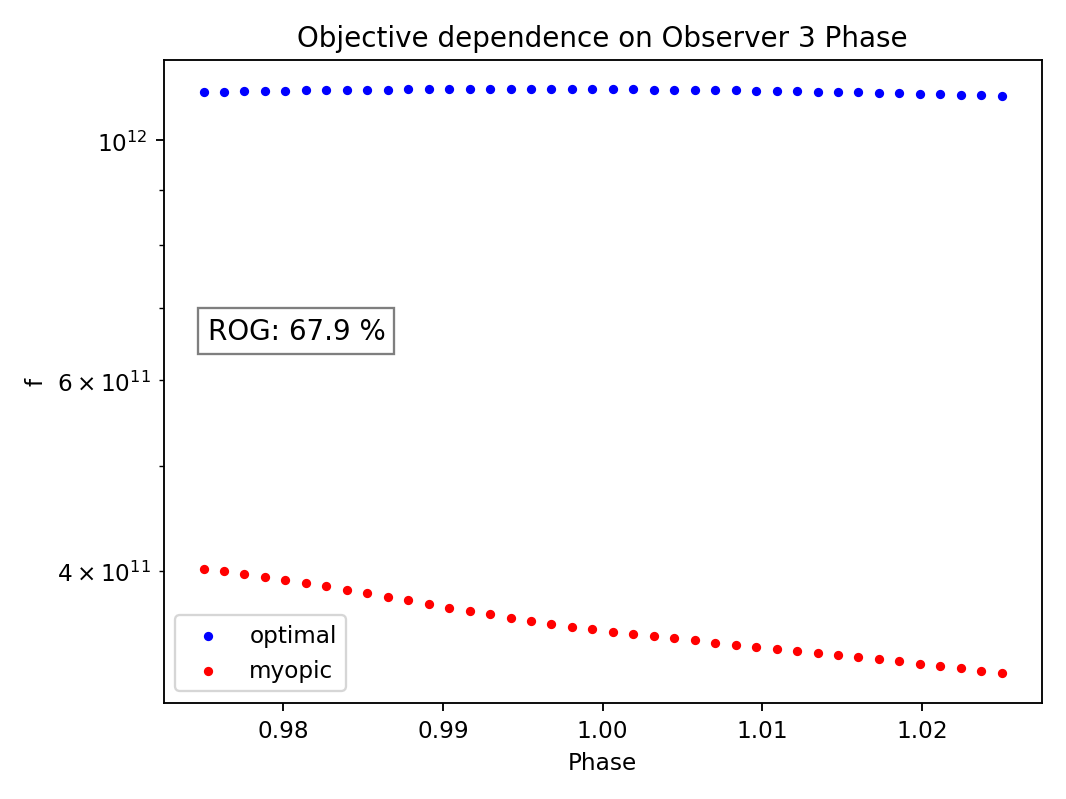}
  \end{minipage}
  \hfill
  \begin{minipage}[b]{0.49\textwidth}
    \includegraphics[width=\textwidth]{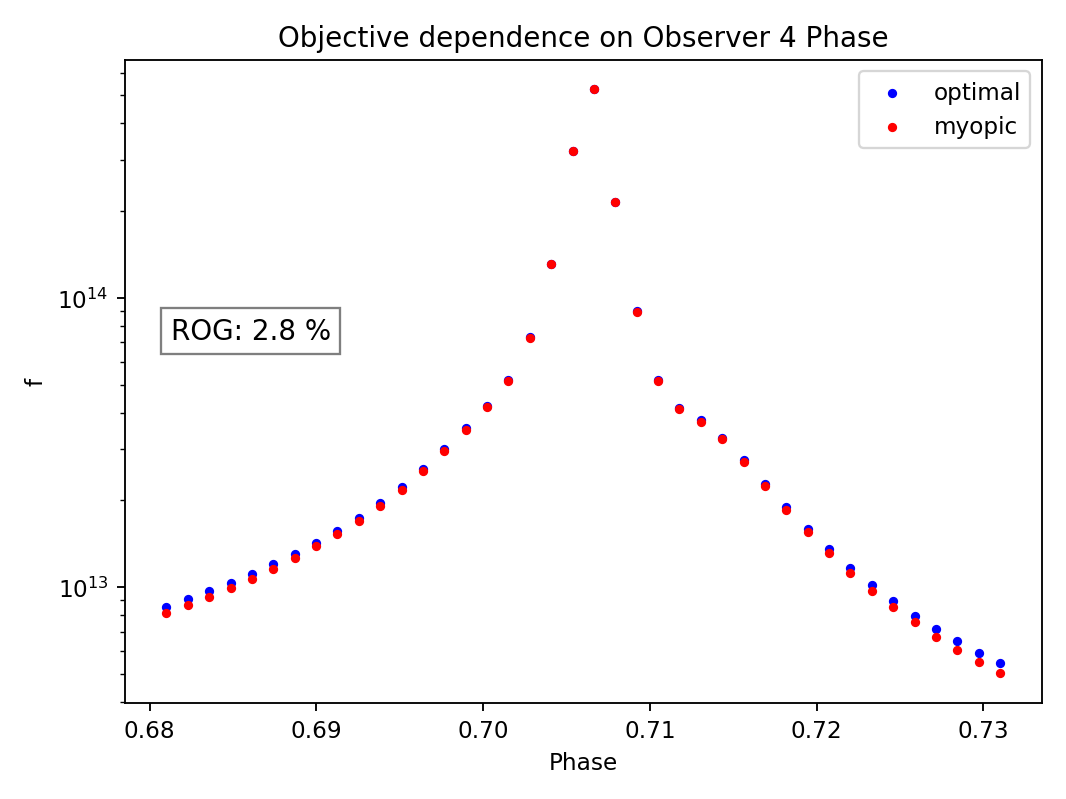}

  \end{minipage}
    \caption{Objective dependence on observer 3 phase (left) and observer 4 phase (right), plotted on a log-linear scale. In each plot, the phase of all other observers are kept constant at 0. The myopic policy is always sub-optimal but approaches optimality if an agent has a close approach to a target.}
    \label{fig:5}
\end{figure}

The orbits in Figure \ref{fig:4} were carefully chosen to avoid intersections between agents and targets, since this would lead to exploding coefficients in the lower level linear program in Equation \ref{eq:3}. Observer 4 has a close approach to a target on a L2 Lypaunov orbit while observer 3 does not have any close approaches to targets. Figure \ref{fig:5} sweeps $P = 40$ phase values for observers 3 and 4 near each's optimal phase and shows the resulting objective when using optimal vs. myopic control. See Table \ref{tab:1} for optimal phases. To quantify the difference between the optimal and myopic objectives, we compute the relative gap between them, averaged over $P$, the number of phases considered. 
\begin{align}
    \text{Relative Optimality Gap (ROG)} \triangleq \frac{1}{P} \sum_{i=1}^{P} \frac{f^i_{\text{opt}} - f^i_{\text{myop}}}{ f^i_{\text{opt}} }
\end{align}
Here, $f^i_{\text{opt}}$ and $f^i_{\text{myop}}$ are the objectives calculated using the $i$-th phase and using optimal or myopic control, respectively.

The myopic policy is consistently sub-optimal near observer 3's optimal phase because it does not offer close approaches to any targets. However, the myopic policy is close to optimal for observer 4 which has a close approach to a target. Close approaches to targets are reflected as large information coefficients relative to other timesteps. As a result, optimal control chooses to observe the close target because of the large information coefficients. Myopic control chooses to observer the close target by definition. Because the coefficient dwarfs others, if both myopic and optimal control are chosen such that they observe a target at timestep(s) with a large information coefficient, then there will be a small relative gap between the final objective achieved by both policies.
We conclude that computing the optimal control is justified for non-close approaches as it provides considerable improvement to the objective over myopic control. Further work is required to quantify and delineate a non-close and close approach.

\subsection{Case Study 3: Using the MaxMin Objective to Diversify Control}
\begin{figure}[H]
    
  \centering
  \begin{minipage}[b]{0.49\textwidth}
    \includegraphics[width=\textwidth]{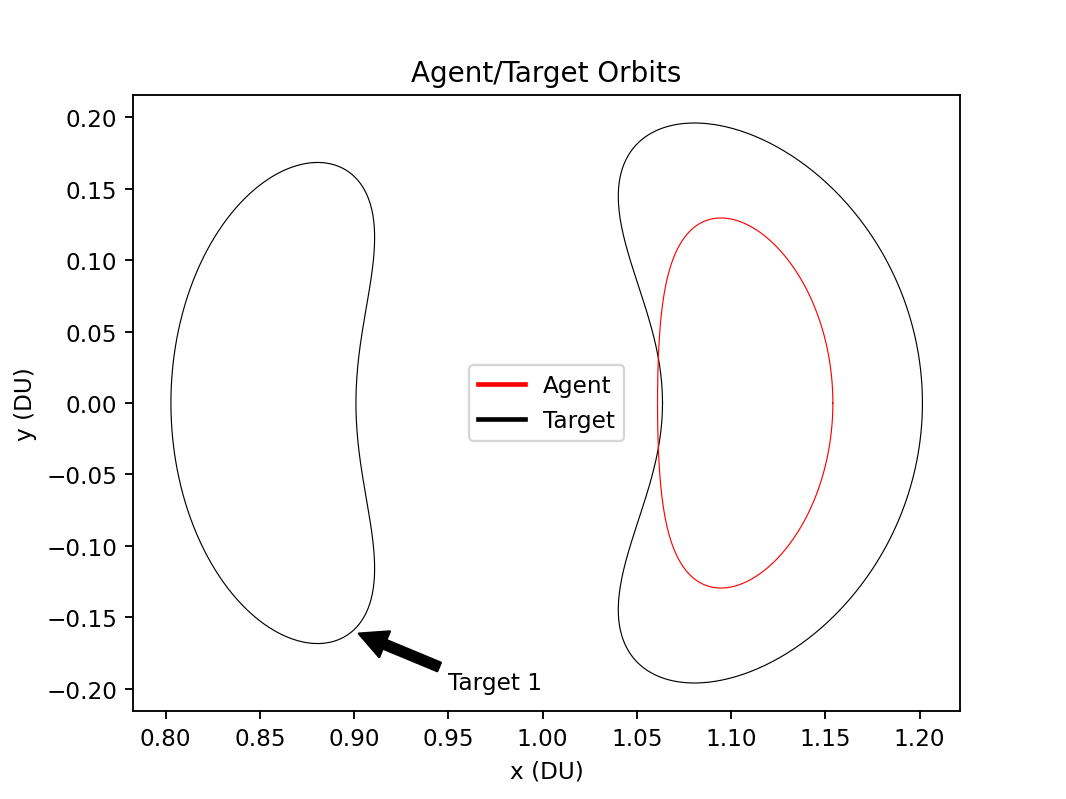}
  \end{minipage}
  \hfill
  \begin{minipage}[b]{0.49\textwidth}
    \includegraphics[width=\textwidth]{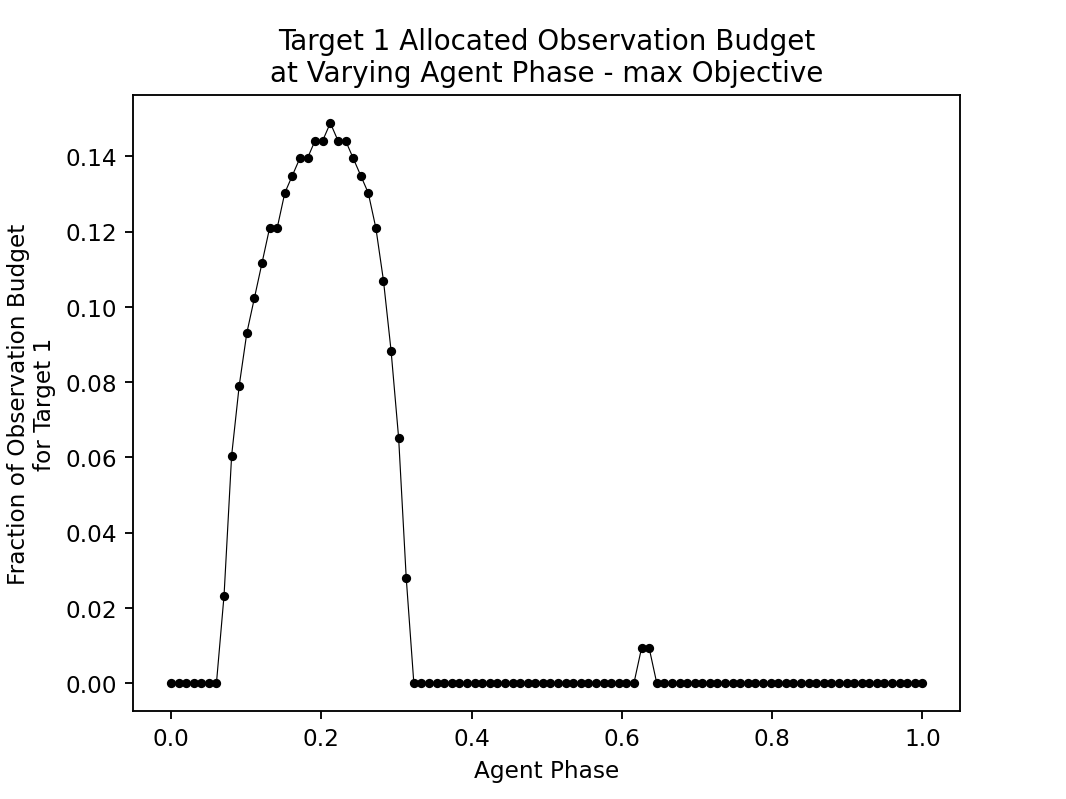}

  \end{minipage}
    \caption{An Agent/Target configuration (left) illustrating how the objective in Equation \ref{eq:1} can be biased. Target 1 is allocated less than 15\% of the observation budget for any phase of the agent.}
    \label{fig:7}
\end{figure}

In the previous sections, the objective was to maximize the total cumulative information (i.e. solve Equation \ref{eq:3} using the objective from Equation \ref{eq:1}). We notice that this objective incentivizes discriminatory control, where agents are not explicitly required to collaborate the achieve the optimal objective. Empirically, we find using this objective leads to optimal control that neglects some target(s) for a large fraction of timesteps, focusing a majority on its effort on a select few targets providing more information. This effect is apparent when an agent/target configuration is chosen such that a subset of targets are favorable for observation (i.e. average distance to agent is smaller). See Figure \ref{fig:7} for an illustration of a configuration where this occurs. Two targets are placed orbits that are closer/further from the agent orbit. The resulting control is biased towards the closer target at all phases of the agent. Note that at the optimal phase $\boldsymbol{x}^* = [0.54]$, the neglected target gets zero observation time.

\begin{figure}[H]
\centering\includegraphics[width=0.9\textwidth]{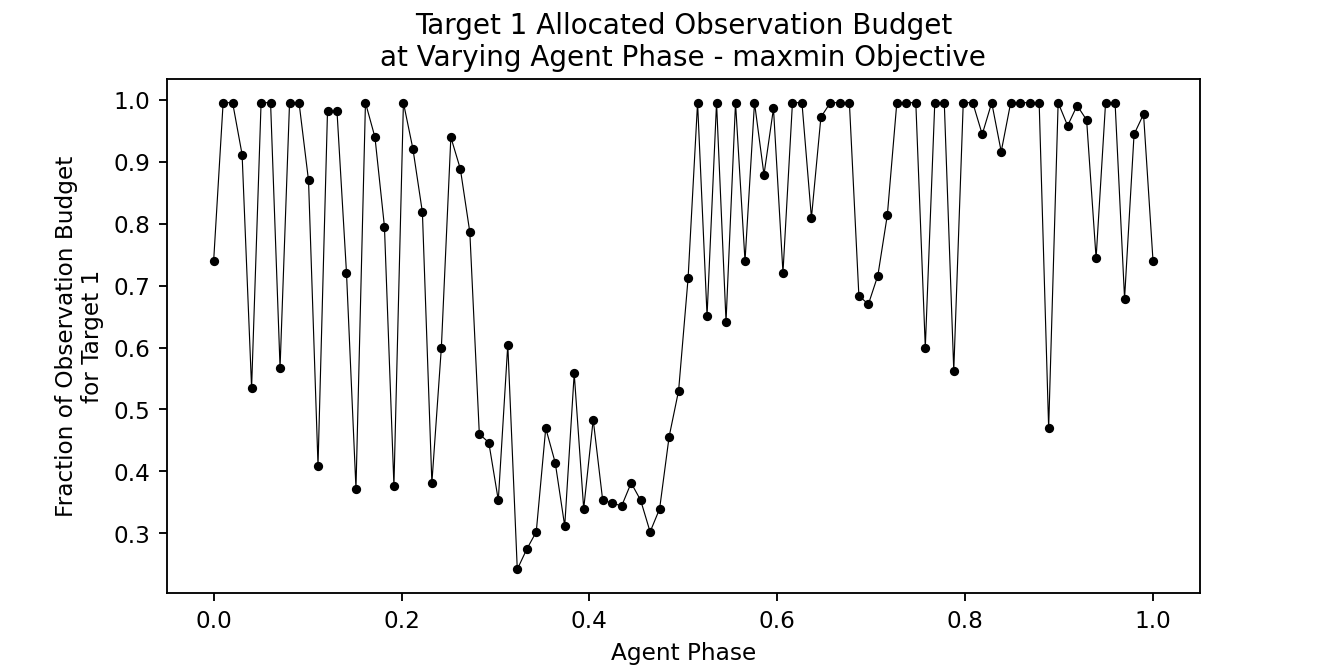}
	\caption{Target 1 observation budget using the maxmin objective in Equation \ref{eq:2}. Target 1 is allocated at least 30\% of the observation budget across all phases.}
	\label{fig:8}
\end{figure}

To reduce this bias against target 1, we employ the maxmin objective in Equation \ref{eq:2}. The resulting allocation budget at various phases is shown in Figure \ref{fig:8}. At the optimal phase $\boldsymbol{x}^* = [0.42]$, target 1 now receives greater than 30\% of the observation budget.

With the maxmin objective, the bias is now towards target 1, although not as severe as the bias towards target 2 with the max objective from Equation \ref{eq:1}. 

\section{Conclusion}
In this paper we present a method to concurrently optimize sensor phasing and tasking in a constellation. Given a set of orbits and number of requested agents, the program returns the phase vector and tasking schedule that maximizes one of the objectives in Equation \ref{eq:1} or Equation \ref{eq:2}. We find that the upper-level objective is smooth-enough for gradient-based techniques. Further work is required to quantify the performance gap between gradient-based and derivative-free optimization techniques. Furthermore, we show that the greedy algorithm for the objective in Equation \ref{eq:1} returns the optimal solution. We show that a myopic policy can perform close to optimal in some cases. The optimality gap is small for constellations with close approaches between targets and observers and the gap grows as this distance grows. Further work is required to quantify close-approaches and their effect on the optimality gap.

\newpage
\bibliographystyle{AAS_publication}   
\bibliography{references}   

\section{Acknowledgement}
This research was supported by the Air Force Office of Scientific Research (AFOSR), as part of the Space University Research Initiative (SURI), grant FA9550-22-1-0092 (grant principal investigator: J. Crassidis from University at Buffalo, The State University of New York). The authors gratefully acknowledge these supporters in the pursuit of this research.

\end{document}